\documentclass{elsart}               

\usepackage{amsmath}
\usepackage{amssymb}
\newtheorem{proposition}{Proposition}[section]

\usepackage{graphicx}          
\usepackage{xcolor}
\usepackage{url}

\begin{document}

\begin{frontmatter}

\title{An Efficient Approximation of the Kalman Filter for Multiple Systems Coupled via Low-Dimensional Stochastic Input} 
\author[Princeton]{Leonid Pogorelyuk}\ead{leonidp@princeton.edu},    
\author[Princeton]{Clarence W. Rowley},               
\author[Princeton]{N. Jeremy Kasdin}  

\address[Princeton]{Department of Mechanical and Aerospace Engineering, Princeton University}  


\begin{abstract}                          
We formulate a recursive estimation problem for multiple dynamical systems
coupled through a low dimensional stochastic input, and we propose an
efficient sub-optimal solution. The suggested approach is an approximation
of the Kalman filter that discards the off diagonal entries of
the correlation matrix in its
``update'' step. The time complexity associated
with propagating this approximate block-diagonal covariance is linear
in the number of systems, compared to the cubic complexity of the
full Kalman filter. The stability of the proposed block-diagonal filter
and its behavior for a large number of systems are analyzed in some
simple cases. It is then examined in the context of electric field
estimation in a high-contrast space coronagraph, for which it was
designed. The numerical simulations provide encouraging results for
the cost-efficiency of the newly suggested filter.
\end{abstract}

\end{frontmatter}

\section{Introduction}

The complexity of the Kalman filter (KF)\cite{kalman1960new} is cubic in the size of the state estimate of the  underlying system. It therefore becomes infeasible for large dimensional system such as atmospheric models\cite{todling1994suboptimal} or when the computational resources are limited as in a space telescopes\cite{pogorelyuk2019dark}.

Various assumptions on the dynamical system have been introduced to reduce the
computational load of the KF. Reduced-order
filters\cite{nagpal1987reduced,cane1996mapping} maintain a low-rank covariance
matrix which can be efficiently approximated in real time\cite{evensen1994sequential}. These methods propagate the state along carefully chosen but low-dimensional subspaces\cite{cohn1995approximate,lermusiaux1999data}, suitable for discrete approximations of spatially continuous physical systems. Another common assumption, for intrinsically discrete systems such as neural networks\cite{puskorius1991decoupled}, robot localization\cite{glielmo1999interlaced} and multiple sensors with dynamical states\cite{sun2004multi}, is that the process noise of the sub-systems is uncorrelated. 

This paper, however, focuses on an estimation problem that has both elements: a continuous physical system (low-order variations of an electric field in a coronagraph) and a high-dimensional discrete state (camera pixels) with independent dynamics and measurement noise. In section~\ref{sec:filter_definition} we introduce and analyze a computationally  efficient sub-optimal approximation of the KF. It discards the off-diagonal blocks of the covariance matrix during the ``update'' step, thus allowing for high-rank state updates while accounting for low-dimensional cross-correlation between numerous sub-systems. In section \ref{sec:examples} we provide numerical examples and compare the proposed block-diagonal filter to the full KF and a banded filter\cite{parrish1985kalman} (which discards the off-diagonal blocks of the covariance at the ``predict'' step), in the context of high-contrast imaging.

\section{\label{sec:filter_definition}State Estimation for Distributed Systems
with Cross-Correlated Process Noise}

\subsection{Problem Formulation}

Consider a set of sub-systems numbered $i=1,...,n$, each
with its own linear dynamics, $F^{(i)}$, and independent process
and measurement noises, $v^{(i)}$ and $w^{(i)}$ respectively. In
addition, the sub-systems are coupled via a stochastic input, $u$,
so that
\begin{alignat}{1}
x_{k+1}^{(i)} &= F_{k}^{(i)}x_{k}^{(i)}+v_{k}^{(i)}+G^{(i)}u_{k}\label{eq:state_equation}\\
y_{k}^{(i)} &= H_{k}^{(i)}x_{k}^{(i)}+w_{k}^{(i)}\label{eq:measurement_equation}
\end{alignat}
where $x_{k}^{(i)},v_{k}^{(i)}\in\mathbb{R}^{c}$, $y_{k}^{(i)},w_{k}^{(i)}\in\mathbb{R}^{d}$,
$u_{k}\in\mathbb{R}^{r}$, and the matrices $F_{k}^{(i)},H_{k}^{(i)},G^{(i)}$
have suitable dimensions. The process and measurement noise are white in time, mutually independent, and normally distributed:
\begin{equation*}
  \label{eq:5}
  w_{k}^{(i)}\sim\mathcal{N}\big(0,R_{k}^{(i)}\big),\:
  v_{k}^{(i)}\sim\mathcal{N}\big(0,V_{k}^{(i)}\big),\:
  u_{k}\sim\mathcal{N}\big(0,U_{k}\big).
\end{equation*}

The full system, obtained by stacking the corresponding vectors of all
sub-systems ($x_{k}=\left[x_{k}^{(i)}\right]_{i=1}^{n}\in\mathbb{R}^{cn}$, ...),
is written
\begin{alignat*}{1}
x_{k+1} &=  F_{k}x_{k}+v_{k}+Gu_{k}\\
y_{k} &=  H_{k}x_{k}+w_{k}.
\end{alignat*}
Here $F_k$ and $H_k$ (as well as $V_k$ and $R_k$) are block diagonal with the
corresponding sub-system matrices ($F_k^{(i)}$, etc.) on the main diagonal, and
$G$ is composed by stacking all of $G^{(i)}$. Assuming $c,d,r\ll n$, the
computational complexity of a KF corresponding to this system is
$\mathcal{O}(n^3)$ despite all the above matrices being either sparse or
low-rank. The main idea of this paper, as discussed in the next
  section, is to exploit a block-diagonal approximation of the covariance in
the ``update'' step of the Kalman filter.  This approximation allows propagating the state estimate in $\mathcal{O}(n)$ complexity while partially accounting for the cross-correlation due to the coupling~$u_k$.

\subsection{\label{sub:diagonal_filter}A Computationally Efficient Block-Diagonal
Filter}

Formally, the block-diagonal approximation is defined via a projection
operator $\mathcal{D}$ over the space of $nc\times nc$ matrices,
\begin{equation*}
\mathcal{D}\left\{ \begin{bmatrix}M_{1,1} & \cdots & M_{1,i} & \cdots\\
\vdots & \ddots & \vdots\\
M_{i,1} & \cdots & M_{i,i} & \cdots\\
\vdots &  & \vdots & \ddots
\end{bmatrix}\right\} =\begin{bmatrix}M_{1,1} & \cdots & 0 & \cdots\\
\vdots & \ddots & \vdots\\
0 & \cdots & M_{i,i} & \cdots\\
\vdots &  & \vdots & \ddots
\end{bmatrix},
\end{equation*}
where $M_{i,i}\in\mathbb{R}^{c\times c}$ . Note that for any symmetric
positive-definite matrix $M$, $\mathcal{D}\left\{ M\right\} $ is also
symmetric positive-definite.

The proposed filter relies on a block-diagonal approximations of the ``update''-step covariance matrices, $P_{k|k}$, denoted by $\tilde{P}_{k|k}=\mathcal{D}\big\{ \tilde{P}_{k|k}\big\}$. Specifically, the projection occurs at the ``update'' step of the KF,
\begin{alignat}{1}
\tilde{P}_{k|k-1} &= F_{k}\tilde{P}_{k-1|k-1}F_{k}^{T}+V_{k}+GU_{k}G^{T}\label{eq:P_tilde_predict}\\
\tilde{S}_{k} &= H_{k}\tilde{P}_{k|k-1}H_{k}^{T}+R_{k}\\
\tilde{K}_{k} &= \tilde{P}_{k|k-1}H_{k}^{T}\tilde{S}_{k}^{-1}\label{eq:new_gain}\\
\tilde{P}_{k|k} &= \mathcal{D}\left\{ \left(I-\tilde{K}_{k}H_{k}\right)\tilde{P}_{k|k-1}\right\} \label{eq:P_tilde_update}\\
\tilde{x}_{k} &= F_{k}\tilde{x}_{k-1}+\tilde{K}_{k}\left(y_{k}-H_{k}F_{k}\tilde{x}_{k-1}\right).\label{eq:x_tilde_update}
\end{alignat}
Note that this filter does not keep track of its true error covariance, which
is neither $P_{k|k}$, the covariance of the optimal filter, nor its
approximation, $\tilde{P}_{k|k}$. However, it achieves a significantly faster performance without discarding the cross-correlation altogether (as would occur if the projection was performed at eq.~(\ref{eq:P_tilde_predict}), resulting in a banded filter\cite{parrish1985kalman}).

Indeed, it is possible to advance
the state and covariance estimates (eqs.~(\ref{eq:P_tilde_update})
and (\ref{eq:x_tilde_update})) in $\mathcal{O}\left(nr^{2}\right)$
time complexity (for a code example, see \cite{pogorelyuk2019block}):
\begin{equation}
\tilde{P}_{k|k}=A_{k}+\mathcal{D}\left\{ B_{k}C_{k,1}B_{k}^{T}+GC_{k,2}G^{T}+GC_{k,3}B_{k}^{T}+B_{k}C_{k,3}G^{T}\right\} ,\label{eq:efficient_representation}
\end{equation}
where
\begin{alignat}{2}
A_{k} &= L_{k}-L_{k}H_{k}^{T}M_{k}^{-1}H_{k}L_{k} && \in\mathbb{R}^{cn\times cn}\\
B_{k} &= L_{k}H_{k}^{T}M_{k}^{-1}H_{k}G && \in\mathbb{R}^{cn\times r}\\
C_{k,1} &= \left(U_{k}^{-1}+N_{k}\right)^{-1} && \in\mathbb{R}^{r\times r}\\
C_{k,2} &= U_{k}\left(N_{k}C_{k,1}N_{k}-N_{k}\right)U_{k}+U_{k} \qquad&& \in\mathbb{R}^{r\times r}\\
C_{k,3} &= U_{k}N_{k}C_{k,1}-U_{k}, && \in\mathbb{R}^{r\times r}
\end{alignat}
and
\begin{alignat}{2}
L_{k} &= F_{k}\tilde{P}_{k-1|k-1}F_{k}^{T}+V_{k} \qquad && \in\mathbb{R}^{cn\times cn}\label{eq:L_def}\\
M_{k} &= H_{k}L_{k}H_{k}^{T}+R_{k} && \in\mathbb{R}^{dn\times dn}\label{eq:M_def}\\
N_{k} &= G^{T}H_{k}^{T}M_{k}^{-1}H_{k}G. && \in\mathbb{R}^{r\times r}\label{eq:N_def}
\end{alignat}
Here the matrices $A_{k},L_{k},M_{k}$ are all block-diagonal and computed
in $\mathcal{O}(n)$. $N_{k}$ and $B_{k}$ can be computed in $\mathcal{O}(nr^{2})$
(performing the sparse multiplications first), and $C_{k,1},C_{k,2},C_{k,3}$
are then all computed in $\mathcal{O}(r^{3})$ via operations on $r\times r$
matrices. Furthermore, the $\mathcal{O}(n)$ non-zero elements of $\tilde{P}_{k|k}$
are computed in $\mathcal{O}(r^{2})$ each, by performing multiplications
of the corresponding $r$ dimensional rows and columns of matrices
inside $\left\{ \cdot\right\} $ in eq.~(\ref{eq:efficient_representation}). Eqs.~(\ref{eq:efficient_representation})--(\ref{eq:N_def}) were derived by several applications of the matrix inversion lemma and grouping terms based on whether they are sparse. 

The sub-optimal gain can be expressed as
\begin{equation}
\tilde{K}_{k}=\left(L_{k}-B_{k}C_{k,1}G^{T}-GC_{k,3}G^{T}\right)H_{k}^{T}M_{k}^{-1},\label{eq:new_gain_efficient}
\end{equation}
which allows updating the estimate in eq.~(\ref{eq:x_tilde_update}) in
$\mathcal{O}(nr)$ by only performing vector multiplications. The parameters of the posterior distribution of the coupling input, $u_k \sim \mathcal{N}\left(\mu_{u,k},\Sigma_{u,k}\right)$, can be approximated via,
\begin{alignat*}{1}
\mu_{u,k} &\approx -C_{k,3}G^{T}H_{k}^{T}M_{k}^{-1}\left(y_{k}-H_{k}F_{k}\hat{x}_{k-1}\right),\\
\Sigma_{u,k} &\approx \left(U_{k}^{-1}+N_{k}\right)^{-1}=C_{k,1}.
\end{alignat*}

\subsection{\label{sec:filter_properties}Steady-State Block-Diagonal Filter}
To gain some theoretical insight into the behavior of the sub-optimal block-diagonal filter, we consider the case of fixed matrices ($F_k=F,\:\forall k$, etc.). Although it should be noted that a fixed-gain approach is more suitable in this case, we suspect that our conclusions could be generalized to systems with periodic or random coefficients (similarly
to \cite{bittanti1988difference} and~\cite{de1984optimal}).

We first define the limit as $k \rightarrow \infty$ of the ``predict''-step covariance approximation, $\tilde{P}_{k|k-1}$ in eq.~(\ref{eq:P_tilde_predict}). It behaves similarly to the steady-state covariance of the full KF whose limit is given by the positive-definite solution of Ricatti equation (for any $Q>0$),
\begin{equation}
P_{\text{-}}=F\left(P_{\text{-}}-P_{\text{-}}H^{T}\left(HP_{\text{-}}H^{T}+R\right)^{-1}HP_{\text{-}}\right)F^{T}+Q.\label{eq:ricatti}
\end{equation}

\begin{proposition}
  \label{par:stability_proposition}
If the pair $\left[F,H\right]$ is detectable and $Q>0$, then
$\tilde{P}{}_{k|k-1}$ converges to the unique positive-definite solution
of
\begin{equation}
\tilde{P}_{\text{-}}=F\mathcal{D}\left\{ \tilde{P}_{\text{-}}-\tilde{P}_{\text{-}}H^{T}\left(H\tilde{P}_{\text{-}}H^{T}+R\right)^{-1}H\tilde{P}_{\text{-}}\right\} F^{T}+Q\label{eq:diagonal_ricatti}
\end{equation}
for any $\tilde{P}{}_{0|0}>0$. (For proof see Appendix~\ref{sub:stability_proof}).
\end{proposition}

We may now say that the steady state covariance, $P_{\text{-}}(Q)$, and its block-diagonal approximation, $\tilde{P}_{\text{-}}(Q)$, are functions of the process noise, $Q$. The same goes for their ``update''-step equivalents, $P_{\text{+}}(Q),\tilde{P}_{\text{+}}(Q)$ (note that $\tilde{P}_{\text{+}}(Q)=\mathcal{D}\left\{\tilde{P}_{\text{+}}(Q)\right\}$).

Is the filter in eqs.~(\ref{eq:P_tilde_predict})-(\ref{eq:x_tilde_update}) stable and does $\tilde{P}_{\text{+}}(V+GUG^T)$ become a better approximation of $P_{\text{+}}(V+GUG^T)$ as the number of sub-system, $n$, increases? Unfortunately, the answer to both questions is not necessarily. The next proposition provides a condition that guarantees stability when the coupling input is small, but the limit as $n \rightarrow \infty$ is more subtle.
As shown in section~\ref{sub:decoupling}, the Kalman filter doesn't necessarily become ``decoupled'' even though adding more sub-systems provides more information about the coupling input. We will show that in some (but not all) cases, $u_k$ tends to behave as ``deterministic'' in the limit $n \rightarrow \infty$, i.e. $P_\text{+}(V+GUG^T)$ resembles $P_\text{+}(V)$, although the characterization of these cases is beyond the scope of this paper.

The block-diagonal approximation, $\tilde{P}_{\text{+}}$,  behaves somewhat differently. To quantify its behavior as $n \rightarrow \infty$ we define,
\begin{equation}
C=\left(U^{-1}+G^{T}H^{T}\left(H\left(F\grave{P}_{\text{+}}F^{T}+V\right)H^{T}+R\right)^{-1}HG\right)^{-1}, \label{eq:C_def}
\end{equation}
where $\grave{P}_{\text{+}}=P_{\text{+}}(\mathcal{D}\left\{ V+GUG^{T}\right\} )$ would be the KF covariance if the input $u_k$ was independent between the sub-systems (alternatively, this is the covariance approximation of the banded filter). Note that since the ``old'' elements of $\grave{P}_{\text{+}}$ are not affected by addition of new sub-systems ($\grave{P}_{\text{+}}$ is block-diagonal), the second term on the right-hand-side of eq.~(\ref{eq:C_def}) is non-decreasing with the addition of new systems. Hence, $C$ as a function of the number of systems is non-increasing (if the sub-systems are the all same, $\left\Vert C\right\Vert =\Theta(n^{-1})$).

\begin{proposition}
\label{par:convergence_proposition}
Let $F_{c}(Q)=\left(I-K(Q)H\right)F$ be the (stable) closed loop transition matrix corresponding to the asymptotic gain $K=P_{\text{-}}H^{T}\left(HP_{\text{-}}H^{T}+R\right)^{-1}$, and let the effects of coupling be described via
\begin{alignat*}{1}
\epsilon^{(i)}&=\left\Vert G^{(i)}CG^{(i)T}\right\Vert _{F},\\
\eta&=\left\Vert GUG^{T}\right\Vert _{F}.
\end{alignat*}
\begin{enumerate}
    \item {
    For small enough $\epsilon^{(i)}$, the difference $\Delta \tilde{P}^{(i)}_{\text{+}}$  between the $i$-th block on the diagonals of $P_{\text{+}}(V)$ and $\tilde{P}_{\text{+}}(V+GUG^T)$ satisfies
    \begin{equation*}
        \left\Vert \Delta \tilde{P}^{(i)}_{\text{+}}\right\Vert _{F} = O(\epsilon^{(i)}).
    \end{equation*}
    }
    \item {
    Let $\Delta P _{\text{-}}$ stand for either $P_{\text{-}}(V+GUG^{T})-P_{\text{-}}(V)$ or $\tilde{P}_{\text{-}}(V+GUG^{T})-P_{\text{-}}(V)$ (and similarly, $\Delta F _{c}$). For small enough $\eta$,
    \begin{alignat*}{1}
    \left\Vert \Delta P_{\text{-}}\right\Vert _{F}=&O(\eta),\\
    \left\Vert \Delta F_{c}\right\Vert _{F}=&O(\eta).
    \end{alignat*}
    }
\end{enumerate}
The exact conditions on $\epsilon^{(i)}$, $\eta$ and the bounds on $\left\Vert \Delta \tilde{P}^{(i)}_{\text{+}}\right\Vert _{F}$, $\left\Vert \Delta P_{\text{-}}\right\Vert _{F}$ and $\left\Vert \Delta F_{c}\right\Vert _{F}$ are given in Appendix~\ref{sub:convergence_proof} together with the proof of the proposition.
\end{proposition}

The first part of the proposition implies that \textit{if} the covariance of full KF is well characterized by  its main block-diagonal, \textit{then} the block-diagonal filter is its good approximation. The example in sec.~\ref{sub:decoupling} illustrates this point as well as a case when the full KF doesn't become ``decoupled.'' The second part of the proposition ensures that the block diagonal filter remains close to optimal for small process noise. The stability of $\tilde{F}_c$ then follows from standard eigenvalue perturbation theory\cite{bauer1960norms}.

\section{\label{sec:examples}Numerical Results}

\subsection{\label{sub:decoupling}Asymptotic ``Decoupling'' of the Steady State Filter}

\begin{figure*}
\includegraphics[height=7cm]{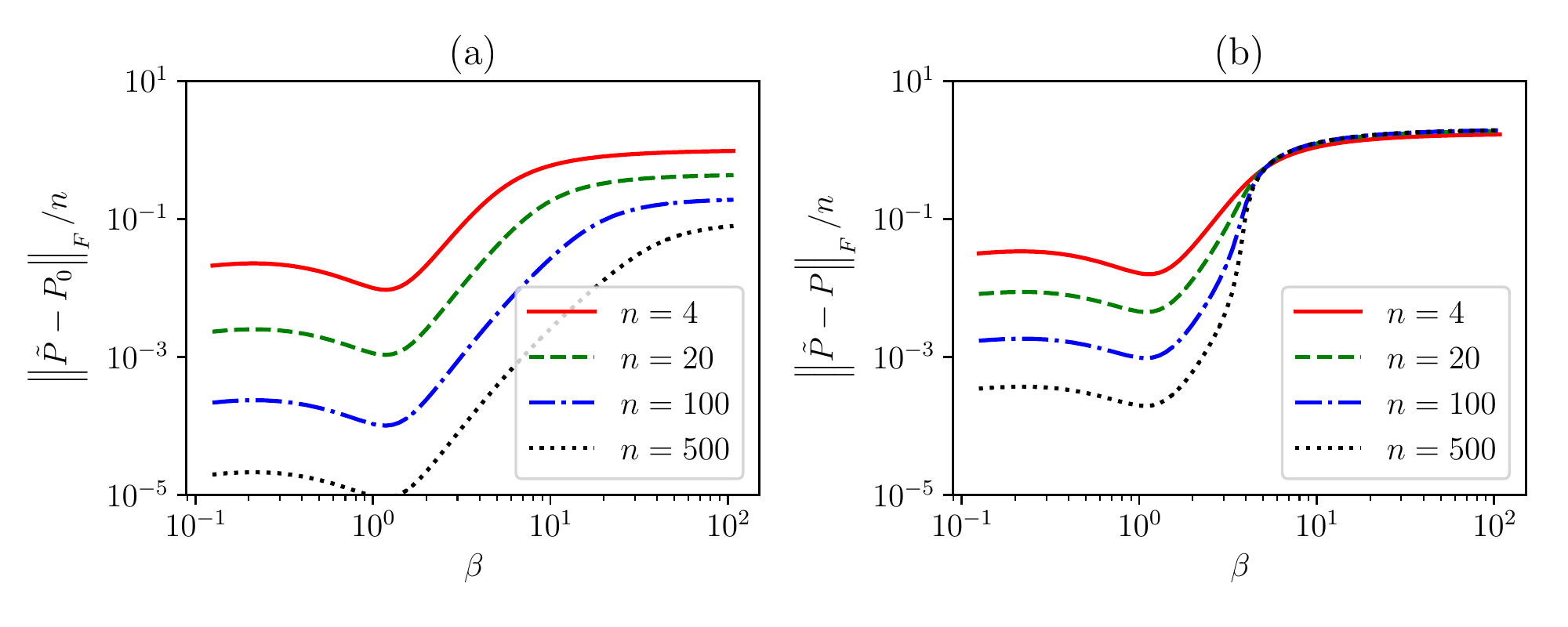}\caption{(a) The difference, scaled by the number of sub-systems, between the approximate covariance of the block-diagonal filter and
the covariance of the ``uncoupled'' KF (for $U=0$), corresponding to the system defined in eqs.~(\ref{eq:example_dynamics})--(\ref{eq:example_coupling}). Each block of the diagonal approximation becomes closer to the corresponding block of the ``uncoupled'' covariance, regardless of the value of the system parameter, $\beta$. (b) The difference between the asymptotic covariance of the full KF and its block-diagonal approximation. Up to a certain critical $\beta $, the full KF becomes increasingly ``decoupled'' and well approximated by the block-diagonal filter.}
\label{fig:convergence_example}

\end{figure*}

Below we illustrate both the cases of a KF remaining ``coupled'' and becoming ``decoupled'' as the number of sub-systems, $n$, increases. It is compared to its block diagonal approximation and a sample code using this systems is given in \cite{pogorelyuk2019block}. Consider a
parameter $\beta$ and a set of $n$ identical systems evolving according
to eqs.~(\ref{eq:state_equation}) and (\ref{eq:measurement_equation})
with

\begin{equation}
F^{(i)}=\begin{bmatrix}0.9 & \beta\\
0 & 0.9
\end{bmatrix},\:H^{(i)}=\begin{bmatrix}1 & 1\end{bmatrix},\:V^{(i)}=\begin{bmatrix}1 & 0\\
0 & 1
\end{bmatrix},\:R^{(i)}=\begin{bmatrix}1\end{bmatrix},\label{eq:example_dynamics}
\end{equation}
and coupling
\begin{equation}
U=\begin{bmatrix}1\end{bmatrix},\:G^{(i)}=\begin{bmatrix}1\\
1
\end{bmatrix}.\label{eq:example_coupling}
\end{equation}

For brevity we denote the steady-state covariance of the Kalman filter (see eq.~(\ref{eq:ricatti}) and discussion below proposition \ref{par:stability_proposition}) as $P=P_{\text{+}}(V+UGU^T)$, its block-diagonal approximation as $\tilde{P}=\tilde{P}_{\text{+}}(V+UGU^T)$ and the covariance of the filter corresponding to the uncoupled case (no process noise) as $P_{0}=P_{\text{+}}(V)=\tilde{P}_{\text{+}}(V)$.

Since both $\tilde{P}$ and $P_{0}$ are block-diagonal, we may compare
them by the average distance between their blocks, $\big\Vert \tilde{P}-P_{0}\big\Vert _{F}/n$.
Consequently, we also consider $\left\Vert \tilde{P}-P\right\Vert _{F}/n$
as a measure of how close the full KF gets to its approximation. Figure
\ref{fig:convergence_example} compares the two measures for various
values of the parameter $\beta$ and exponentially increasing $n$.
According to proposition \ref{par:convergence_proposition}, $\big\Vert \tilde{P}-P_{0}\big\Vert _{F}/n\rightarrow0,\:\forall\beta$, i.e., each block of the covariance approximation, $\tilde{P}$, becomes close to the covariance of the corresponding sub-system of the uncoupled KF (fig.~\ref{fig:convergence_example}(a)).

The normalized difference between $P$ (full KF) and $P_0$ (uncoupled KF) also decreases with $n$ but only for some values of $\beta$ (fig.~\ref{fig:convergence_example}(b)). In these cases the full KF becomes ``decoupled'' and $P$ becomes close to block-diagonal. For other $\beta$, above what appears to be some critical value, $\left\Vert \tilde{P}-P\right\Vert _{F}$ (and $\left\Vert P_{0}-P\right\Vert _{F}$)
grows as $n$, similarly to $\left\Vert GUG^{T}\right\Vert$. Although not shown explicitly, the true covariance of the block-diagonal filter (the covariance of the error of the state estimate obtained by using the sub-optimal gain of the block-diagonal filter), also converges to $P_{0}$ but
only in the region where $\big\Vert \tilde{P}-P\big\Vert _{F}/n\rightarrow0$.

\subsection{\label{sub:speckle_drift}Speckle Drift in High-Contrast Imaging}

To observe planets orbiting stars outside our solar system, modern telescopes employ coronagraphs---devices that block most of the light from an observed star to create an image with high-contrast regions (a dark hole) where faint objects can be detected. Contrasts on the order of $10^8$ between the exo-planet and its host star can only be achieved from space-based telescopes such as the proposed Wide-Field InfraRed Survey Telescope (WFIRST) Coronagraph Instrument (CGI)\cite{demers2015requirements}.
Even in space, the residual star light (speckles) does not remain static over the tens of hours required to collect enough photons to detect a planet, imposing tight requirements on the stability of the optics\cite{shaklan2011stability}.

The authors have proposed actively maintaining the dark hole contrast via real-time estimation of the speckles electric field from intensity measurements in the camera plane\cite{pogorelyuk2019dark}. The method treats the electric field at each pixel of the dark hole as a dynamical system and employs an extended Kalman filter (EKF). In this model, the measurement noise is due to the number of photons detected at each pixel, and the dynamics are a result of stochastic time variation in the speckles field. While the detector noise is uncorrelated between pixels, the speckles in the dark hole are ``smooth'' and can be described by a linear combination of a small number of ``modes''\cite{soummer2012detection}. Below we exploit this ``smoothness'' property to obtain better estimates of the electric field.

For this example we used a simulated model of the Princeton High Contrast Imaging Lab (HCIL) testbed, which can be described as a complex electric field, $x^{(i)}=[\mathrm{Re}E^{(i)},\mathrm{Im}E^{(i)}]^T$, sampled at $n=2088$ pixels. While the initial electric field, $x(0)$, is arbitrary, its increments due to drift are well described by a small number of $r\approx 20$ modes, $F=I$ and $\left\Vert V\right\Vert _{2}\ll\left\Vert U\right\Vert _{2}$. The measurement equation is, however, non-linear, \begin{equation*}
y^{(i)}_{k}\sim\mathrm{Poisson}\left(\left|E^{(i)}_{k}+\Delta E^{(i)}_{k}\right|^{2}\right),
\end{equation*}
where $\Delta E$ are known modulations intended to increase the observability of the corresponding EKF (for its derivation and more details about the system see \cite{pogorelyuk2019dark}). 

\begin{figure}
\begin{center}
\includegraphics[height=7cm]{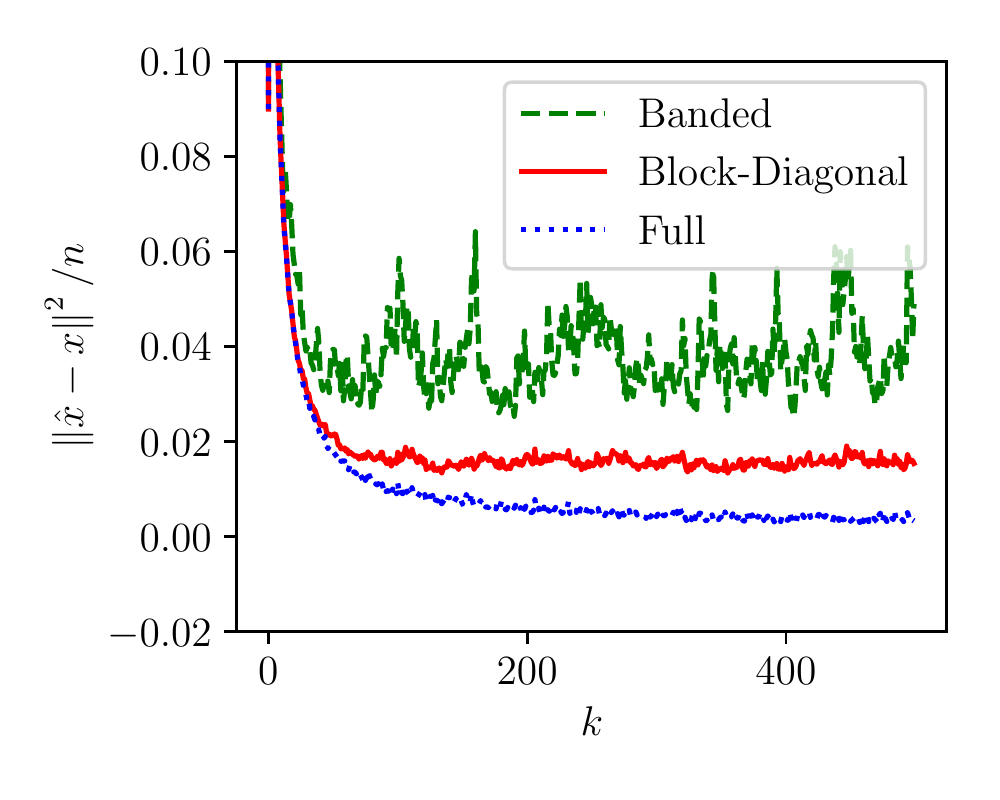}    
\caption{The average squared error of the electric field estimate of the full EKF and of its approximations. The full EKF (dotted blue line) gives the most accurate estimate but is 3 orders of magnitude slower than the block-diagonal filter proposed in this paper (solid red line). The banded EKF\cite{parrish1985kalman,pogorelyuk2019dark} (dashed green line) doesn't exploit the correlation between increments of $x$ and is therefore significantly less accurate but also an order of magnitude faster than the block-diagonal filter.}  
\label{fig:comparison}                                 
\end{center}                                 
\end{figure}

In our simulations, the electric field, $x$, resided in a high dimensional space and the measurements were not correlated between pixels. In \cite{pogorelyuk2019dark} we ignored the smoothness of the electric field increments by propagating an EKF for each pixel independently of the others (a banded EKF\cite{parrish1985kalman}). Figure \ref{fig:comparison} shows that this approach results in a relatively high estimation error. On the other hand, the block-diagonal filter (eqs.~(\ref{eq:efficient_representation})--(\ref{eq:new_gain_efficient})) accounts for the cross-correlation in $x_{k+1}-x_k$, reducing the estimation error by a factor of $2.4$ at the expense of $27$ times longer computation time. The full EKF was further $3.8$ times more accurate but prohibitively $6400$ times slower.

It should be noted that increasing the number of pixels (e.g. by considering multiple wavelengths) would increase the accuracy of both the block-diagonal filter and the full EKF, but would only increase the computation time (per pixel) of the latter. The accuracy of the banded filter doesn't depend on the total number of pixel, but it can be increased by grouping several pixels together into larger sub-systems at the expense of performance. Still, the block-diagonal filter remained more cost-efficient (in terms of accuracy vs. speed) than the banded filter with larger block sizes, although we omit the full analysis here.

\section{Conclusions}

We have formulated an estimation problem for multiple dynamical systems
coupled only through a low dimensional stochastic input. It arises
in the context of electric field estimation in high-contrast telescopes
in the presence of time varying wavefront errors. The computational
requirements of the extended Kalman filter grow as a cube of the number
of sub-systems (pixels in the camera plane), rendering this standard
approach infeasible for an on-board computer of a space telescope.

To reduce the computational burden, we introduced a block-diagonal
approximation of the Kalman filter (see eqs.~(\ref{eq:efficient_representation})--(\ref{eq:new_gain_efficient}) and code example in \cite{pogorelyuk2019block}))
which discards the cross-correlation between the state estimates of each sub-system 
(the off-diagonal blocks of the covariance matrix) after each ``update''
step of the filter. The proposed block-diagonal filter has linear
time complexity and combines the measurements from all sub-systems
during its ``predict'' step.

Our numerical results in sec.~\ref{sub:speckle_drift} suggest
that having an accurate reduced-order model for the drift of the electric field in a coronagraph, can
significantly increase the accuracy of its online estimates.
The most accurate results were obtained by formulating an EKF for
all the dark hole pixels at once, even though this would be computationally infeasible in practice. A banded filter which discards all information regarding
coupling was both the fastest and the least accurate.
The block-diagonal filter, on the other hand, appeared to be the most cost-efficient: it both benefited from the increasing number of sub-systems (similarly to the EKF), and its performance scaled linearly (similarly to the banded filter).

The block-diagonal filter is not guaranteed to be stable even in the
case of time-invariant systems. Although its gain converges to a finite limit (proposition \ref{par:stability_proposition}), the resulting closed-loop transition matrix might turn out to be unstable. Such cases do not seem to be very common and we speculate that the block-diagonal
filter is stable for many types of time-varying systems (e.g. systems with
random zero-mean coefficients).

Finally, one might expect that as the number of sub-systems increases,
the effect of the fixed-dimensional coupling input should become deterministic
(i.e. the covariance of the KF should become close to that of the KF corresponding to the same system with no coupling). Perhaps surprisingly, this is not always the case even for some simple and
``well observable'' systems, as illustrated in sec.~\ref{sub:decoupling}.
Nevertheless, when the KF does become close to ``uncoupled,'' in
the limit of large number of systems, the block-diagonal filter gives
an increasingly accurate approximation of the KF gain (as a consequence
of proposition \ref{par:convergence_proposition}). In that case, the
newly proposed filter yields close-to-optimal state estimates at a fraction
of the computational cost of the Kalman filter.

\begin{ack}                               
This material is based upon work supported by the Army Research Office under award number W911NF-17-1-0512.
\end{ack}

\bibliographystyle{ieeetr}        
\bibliography{block_diagonal}           
                                 
\appendix

\section{Proofs}

\subsection{Proof of proposition~\ref{par:stability_proposition}}
\label{sub:stability_proof}
For the given values of $F$, $H$ and $R$ , define

\begin{equation}
J\left(P,K,W\right)\equiv\left(I-KH\right)\left(FPF^{T}+W\right)\left(I-KH\right)^{T}+KRK^{T}.\label{eq:J_def}
\end{equation}
The Kalman gain $\tilde{K}_{k}$ given by eq.~(\ref{eq:new_gain}) with $Q=V+GUG^T$,
has the optimality property that
\begin{equation}
J\left(\tilde{P}_{k-1|k-1},\tilde{K}_{k},Q\right) =\left(I-\tilde{K}_{k}H\right)\tilde{P}_{k|k-1}\le J\left(\tilde{P}_{k-1|k-1},K,Q\right)\:\forall K,\label{eq:J_optimality}
\end{equation}
and we note that, for any $P^{\#}\le P$ and $W^{\#}\le W$,
\begin{alignat}{1}
J\left(P^{\#},K,W\right)\le & J\left(P,K,W\right),\label{eq:ineq1}\\
J\left(P,K,W^{\#}\right)\le & J\left(P,K,W\right).\label{eq:ineq2}
\end{alignat}

First we assume, without loss of generality, that $\tilde{P}{}_{0|0}$
is block diagonal ($\mathcal{D}\big\{\tilde{P}{}_{0|0}\big\}=\tilde{P}{}_{0|0}$),
and show that $\tilde{P}_{k|k}$ is bounded. Let $\lambda_{max}$
be the largest eigenvalue of $Q$ (hence $\lambda_{max}I\ge Q$).
We define another Kalman filter starting with the same initial $\tilde{P}_{0|0}$,
except that its process noise covariance is $\lambda_{max}I$ instead
of $Q$; we denote the covariance of this Kalman filter by $\Pi_{k|k}$.
This Kalman filter is also a block-diagonal filter, since if $Q$
is block diagonal, then $\mathcal{D}$ acts as the identity in eq.~(\ref{eq:P_tilde_update}).

In order to show that $\tilde{P}_{k|k}$ is bounded, we show by induction
that $\tilde{P}_{k|k}\le\Pi_{k|k}$ for all $k$, starting with
$\Pi_{0|0}=\tilde{P}_{0|0}$. Indeed, if for some $k$, $\tilde{P}_{k-1|k-1}\le\Pi_{k-1|k-1}$
then
\begin{alignat*}{1}
\tilde{P}_{k|k}= & \mathcal{D}\left\{ \underset{K}{\min}J\left(\tilde{P}_{k-1|k-1},K,Q\right)\right\} \le\mathcal{D}\left\{ J\left(\tilde{P}_{k-1|k-1},\kappa_{k},Q\right)\right\} \\
\le & \mathcal{D}\left\{ J\left(\Pi_{k-1|k-1},\kappa_{k},\lambda_{max}I\right)\right\} =\Pi_{k|k}
\end{alignat*}
where $\kappa$ is the optimal gain for $\Pi$ and in the second inequality we have used (\ref{eq:ineq1})
and (\ref{eq:ineq2}). Since the second filter $\Pi_{k|k}$ is
a Kalman filter, it converges (under the assumptions of the proposition\cite{davis1985stochastic}),
and hence $\tilde{P}_{k|k}$ is bounded.

We next show that, if $\tilde{P}_{0|0}=0$, then the sequence $\tilde{P}_{k|k}$
is nondecreasing, and thus converges to a limit (since we have already
shown the sequence is bounded). Let $\left\{ \tilde{P}_{k|k}^{[0]}\right\} $
denote the sequence of $\tilde{P}_{k|k}$ such that $\tilde{P}_{0|0}^{[0]}=0$.
Again proceeding inductively, if for some $k$, $\tilde{P}_{k-1|k-1}^{[0]}\le\tilde{P}_{k|k}^{[0]}$
then
\begin{alignat*}{1}
\tilde{P}_{k|k}^{[0]}= & \mathcal{D}\left\{ \underset{K}{\min}J\left(\tilde{P}_{k-1|k-1}^{[0]},K,Q\right)\right\} \le\mathcal{D}\left\{ J\left(\tilde{P}_{k-1|k-1}^{[0]},\tilde{K}_{k+1}^{[0]},Q\right)\right\} \\
\le & \mathcal{D}\left\{ J\left(\tilde{P}_{k|k}^{[0]},\tilde{K}_{k+1}^{[0]},Q\right)\right\} =\tilde{P}_{k+1|k+1}^{[0]}
\end{alignat*}
where in the second inequality we have used (\ref{eq:ineq1}).
Since $\tilde{P}_{1|1}^{[0]}\ge0=\tilde{P}_{0|0}^{[0]}$, the induction
is anchored, and the sequence $\tilde{P}_{k|k}^{[0]}$\ is non-decreasing
for all $k$. Since this sequence was already shown to be bounded,
it therefore converges to a limit. Let us define
\begin{equation}
\tilde{P}_{\text{-}}\equiv\underset{k\rightarrow\infty}{\lim}\tilde{P}_{k+1|k}^{[0]}=F\left(\underset{k\rightarrow\infty}{\lim}\tilde{P}_{k|k}^{[0]}\right)F^{T}+Q.\label{eq:P_tilde_minus_def}
\end{equation}

We wish to show that $\tilde{P}_{k+1|k}\rightarrow\tilde{P}_{\text{-}}$
for all nonzero $\tilde{P}_{0|0}$. To do this, we first find an upper
bound for $\tilde{P}_{k+1|k}$. Note that, for $\tilde{K}=\tilde{P}_{\text{-}}H^{T}\left(H\tilde{P}_{\text{-}}H^{T}+R\right)^{-1}$,
\begin{alignat}{1}
\tilde{P}_{\text{-}}= & F\mathcal{D}\left\{ \left(I-\tilde{K}H\right)\tilde{P}_{\text{-}}\left(I-\tilde{K}H\right)^{T}+\tilde{K}R\tilde{K}^{T}\right\} F^{T}+Q\label{eq:diagonal_ricatti_K_form}\\
> & F\mathcal{D}\left\{ \left(I-\tilde{K}H\right)\tilde{P}_{\text{-}}\left(I-\tilde{K}H\right)^{T}\right\} F^{T}\nonumber 
\end{alignat}
and therefore there exists $0<\alpha<1$ such that
\[
F\mathcal{D}\left\{ \left(I-\tilde{K}H\right)\tilde{P}_{\text{-}}\left(I-\tilde{K}H\right)^{T}\right\} F^{T}<\alpha\tilde{P}_{\text{-}}.
\]
Moreover, for each $\tilde{P}_{0|0}\ge0$, there exists $\beta>0$
such that $\tilde{P}_{1|0}<\left(\beta+1\right)\tilde{P}_{\text{-}}$.
We now show inductively that $\tilde{P}_{k+1|k}<(\beta\alpha^{k}+1)\tilde{P}_{\text{-}}$
for all $k$. Indeed, if $\tilde{P}_{k|k-1}<\left(\beta\alpha^{k-1}+1\right)\tilde{P}_{\text{-}}$,
then
\begin{alignat*}{1}
\tilde{P}_{k+1|k}= & F\mathcal{D}\left\{ J\left(\tilde{P}_{k-1|k-1},\tilde{K}_{k},Q\right)\right\} F^{T}+Q\\
\le & F\mathcal{D}\left\{ J\left(\tilde{P}_{k-1|k-1},\tilde{K},Q\right)\right\} F^{T}+Q\\
= & F\mathcal{D}\left\{ \left(I-\tilde{K}H\right)\tilde{P}_{k|k-1}\left(I-\tilde{K}H\right)^{T}+\tilde{K}R\tilde{K}^{T}\right\} F^{T}+Q\\
< & F\mathcal{D}\left\{ \left(I-\tilde{K}H\right)\left(\beta\alpha^{k-1}+1\right)\tilde{P}_{\text{-}}\left(I-\tilde{K}H\right)^{T}+\tilde{K}R\tilde{K}^{T}\right\} F^{T}+Q\\
< & \beta\alpha^{k}\tilde{P}_{\text{-}}+\tilde{P}_{\text{-}}
\end{alignat*}
 Therefore $\tilde{P}_{k+1|k}<\left(\beta\alpha^{k}+1\right)\tilde{P}_{\text{-}}$
for all $k$, and $\underset{k\rightarrow\infty}{\lim}\sup\tilde{P}_{k+1|k}=\tilde{P}_{\text{-}}$.

We next find a lower bound for $\tilde{P}_{k+1|k}$, again proceeding
inductively. If for some $k$, $\tilde{P}_{k|k}^{[0]}\le\tilde{P}_{k|k}$,
then using (\ref{eq:ineq2}) and the optimality of the Kalman
gain, we have
\begin{alignat*}{1}
\tilde{P}_{k+1|k+1}^{[0]}= & \mathcal{D}\left\{ J\left(\tilde{P}_{k|k}^{[0]},\tilde{K}_{k+1}^{[0]},Q\right)\right\} \le\mathcal{D}\left\{ J\left(\tilde{P}_{k|k}^{[0]},\tilde{K}_{k+1},Q\right)\right\} \\
\le & \mathcal{D}\left\{ J\left(\tilde{P}_{k|k},\tilde{K}_{k+1},Q\right)\right\} =\tilde{P}_{k+1|k+1}.
\end{alignat*}
Hence, $\tilde{P}_{k|k}^{[0]}\le\tilde{P}_{k|k}$ for all $k$ and
\[
\underset{k\rightarrow\infty}{\lim}\tilde{P}_{k+1|k}=\underset{k\rightarrow\infty}{\lim}F\tilde{P}_{k|k}F^{T}+Q\ge\tilde{P}_{\text{-}},
\]
where we have used eq. (\ref{eq:P_tilde_minus_def}). This concludes
the proof that $\tilde{P}_{k+1|k}\rightarrow\tilde{P}_{\text{-}}$,
for all $\tilde{P}_{0|0}$. Eq. (\ref{eq:diagonal_ricatti}) then follows from eq. (\ref{eq:diagonal_ricatti_K_form}),
and the uniqueness of its solution follows from the uniqueness of
the limit of $\tilde{P}_{k+1|k}$.

\subsection{Proposition~\ref{par:convergence_proposition} and proof}
\label{sub:convergence_proof}

Let $F_{c}(Q)=\left(I-K(Q)H\right)F$ be the (stable) closed loop transition matrix corresponding to the asymptotic gain $K=P_{\text{-}}H^{T}\left(HP_{\text{-}}H^{T}+R\right)^{-1}$ , the following constants be defined for each sub-system as
\begin{alignat*}{1}
\alpha_{1}^{(i)}&=\frac{1}{1-\max\left|\lambda\left\{ F_{c}^{(i)}(V)\right\} \right|^{2}},\\
\alpha_{2}^{(i)}&=\left\Vert F_{c}^{(i)}(V)\right\Vert _{2},\\
\alpha_{3}^{(i)}&=\left\Vert \left(I+P_{\text{-}}^{(i)}(V)H^{(i)T}\left(R^{(i)}\right)^{-1}H^{(i)}\right)^{-1}\right\Vert _{2},\\
\alpha_{4}^{(i)}&=\left\Vert H^{(i)T}\left(H^{(i)}P_{\text{-}}^{(i)}(V)H^{(i)T}+R^{(i)}\right)^{-1}H^{(i)}F^{(i)}\right\Vert _{2},\\
\alpha_{5}^{(i)}&=\left\Vert H^{(i)T}\left(H^{(i)}P_{\text{-}}^{(i)}(V)H^{(i)T}+R^{(i)}\right)^{-1}H^{(i)}\right\Vert _{2},
\end{alignat*}
and the effects of coupling described via
\begin{alignat*}{1}
\epsilon^{(i)}&=\left\Vert G^{(i)}CG^{(i)T}\right\Vert _{F},\\
\eta&=\left\Vert GUG^{T}\right\Vert _{F}.
\end{alignat*}
\begin{enumerate}
    \item {
    Let $\beta^{(i)}$ be a Bauer-Fike bound\cite{bauer1960norms} on the eigenvalue perturbation of $F_{c}^{(i)}(V)$,
    \begin{equation*}
        \left|\lambda\left\{ F_{c}^{(i)}\right\} -\lambda\left\{ F_{c}^{(i)}+\Delta F \right\} \right|\le\beta^{(i)}\left\Vert \Delta F\right\Vert,\:\forall \Delta F 
    \end{equation*}
    If $\alpha_{1}^{(i)}\left(2\alpha_{2}^{(i)}+\alpha_{3}^{(i)}\alpha_{4}^{(i)}\epsilon^{(i)}\right)\alpha_{3}^{(i)}\alpha_{4}^{(i)}\epsilon^{(i)}<1$ and $\max\left|\lambda\left\{ F_{c}^{(i)}(V)\right\} \right|+\beta^{(i)}\alpha_{3}^{(i)}\alpha_{4}^{(i)}\epsilon^{(i)}<1$, then the difference $\Delta \tilde{P}^{(i)}_{\text{+}}$  between the $i$-th block on the diagonals of $P_{\text{+}}(V)$ and $\tilde{P}_{\text{+}}(V+GUG^T)$ satisfies,
    \begin{equation*}
        \left\Vert \Delta \tilde{P}^{(i)}_{\text{+}}\right\Vert _{F}\le\frac{\alpha_{1}^{(i)}\left(\alpha_{2}^{(i)}\right)^{2}\epsilon^{(i)}}{1-\alpha_{1}^{(i)}\left(2\alpha_{2}^{(i)}+\alpha_{3}^{(i)}\alpha_{4}^{(i)}\epsilon^{(i)}\right)\alpha_{3}^{(i)}\alpha_{4}^{(i)}\epsilon^{(i)}}.
    \end{equation*}
    }
    \item {
    Let $\alpha_{j}=\underset{i}{\max}\alpha_{j}^{(i)}$, and $\Delta P _{\text{-}}$ stand for either $P_{\text{-}}(V+GUG^{T})-P_{\text{-}}(V)$ or $\tilde{P}_{\text{-}}(V+GUG^{T})-P_{\text{-}}(V)$ (and similarly, $\Delta F _{c}$). If, $4\alpha_{1}^{2}\alpha_{2}^{2}\alpha_{5}\eta<1$ then,
    \begin{alignat*}{1}
    \left\Vert \Delta P_{\text{-}}\right\Vert _{F}\le&\frac{1-\sqrt{1-4\alpha_{1}^{2}\alpha_{2}^{2}\alpha_{5}\eta}}{2\alpha_{1}\alpha_{2}^{2}\alpha_{5}}\le2\alpha_{1}\eta,\\
    \left\Vert \Delta F_{c}\right\Vert _{F}\le&\alpha_{3}\alpha_{4}\left\Vert \Delta P_{\text{-}}\right\Vert _{F}.
    \end{alignat*}
    }
\end{enumerate}

\subsubsection*{Proof of part (1)}
Throughout the proof we will use shorthands for the steady-state covariance approximation and transition matrix of the block-diagonal filter,
\begin{alignat*}{1}
\tilde{P}_{\text{\ensuremath{\pm}}}= & \tilde{P}_{\text{\ensuremath{\pm}}}\left(V+GUG^{T}\right),\\
\tilde{F}_{c}= & \tilde{F}_{c}\left(V+GUG^{T}\right)
\end{alignat*}
which will be compared to their \emph{decoupled} (no cross-correlated noise terms) counterparts,
\begin{alignat*}{1}
P_{\text{\ensuremath{\pm}}}= & P_{\text{\ensuremath{\pm}}}\left(V\right),\\
F_{c}= & F_{c}(V).
\end{alignat*}

The steady-state estimate of the covariance of the banded filter (which accounts for the coupling noise but ignores the resulting cross correlation) is given after the definition of $C$ in eq.~(\ref{eq:C_def}) as
\begin{equation*}
\grave{P}_{\text{+}}=P_{\text{+}}\left(\mathcal{D}\left\{ V+GUG^{T}\right\} \right).    
\end{equation*}
The sub-optimality of the banded filter implies $\tilde{P}_{\text{+}}\le\grave{P}_{\text{+}}$
(see previous proof). However $\grave{P}_{\text{+}}$ is not a useful
abound on $\tilde{P}_{\text{+}}$, because it doesn't get ``closer''
to $P_{\text{+}}$ since there is no information exchange regarding
$u$ between the filters of each sub-system.

We will now define another ``intermediate'' filter, with its own
block-diagonal approximation of the steady-state covariance matrix,
$\tilde{\Pi}_{\text{+}}$. This filter also disregards the cross-correlation,
except when computing the state estimate of one of the sub-systems
($i=1$, without loss of generality). It will give a bound
on the steady-state covariance approximation of just that one sub-system,
\begin{alignat*}{1}
\tilde{P}_{\text{+}}^{(i)}\le\tilde{\Pi}_{\text{+}}^{(i)}\le\grave{P}_{\text{+}}^{(i)},\: & i=1,\\
\tilde{P}_{\text{+}}^{(i)}\le\tilde{\Pi}_{\text{+}}^{(i)}=\grave{P}_{\text{+}}^{(i)},\: & i\ne1.
\end{alignat*}

The intermediate filter uses a sub-optimal gain,
\[
\tilde{\kappa}=\begin{bmatrix}K_{mix}^{(1,1)} & K_{mix}^{(1,2)} & K_{mix}^{(1,3)} & \cdots\\
 & \grave{K}^{(2)}\\
 &  & \grave{K}^{(3)}\\
 &  &  & \ddots
\end{bmatrix},
\]
defined in two parts: systems with $i\ne1$ use estimates based on
the asymptotic banded filter gain $\grave{K}=K\left(\mathcal{D}\left\{ V+GUG^{T}\right\} \right)$;
the first system ($i=1$) uses the gain
\[
K_{mix}=\left(FP_{mix}F^{T}+V+GUG^{T}\right)H^{T}\left(H\left(FP_{mix}F^{T}+V+GUG^{T}\right)H^{T}+R\right)^{-1},
\]
which does take cross-correlation into account with
\[
P_{mix}=\begin{bmatrix}P_{\text{+}}^{(1)}\\
 & \grave{P}_{\text{+}}^{(2)}\\
 &  & \grave{P}_{\text{+}}^{(3)}\\
 &  &  & \ddots
\end{bmatrix}.
\]

Using some algebraic manipulations, one can show that
\begin{alignat*}{1}
K_{mix}= & \left(FP_{mix}F^{T}+V+\left(I+\left(FP_{mix}F^{T}+V\right)H^{T}R^{-1}H\right)^{-1}GC_{mix}G^{T}\right)\cdot\\
 & \cdot H^{T}\left(H\left(FP_{mix}F^{T}+V\right)H^{T}+R\right)^{-1}
\end{alignat*}
where
\[
C_{mix}=\left(U^{-1}+G^{T}H^{T}\left(H\left(FP_{mix}F^{T}+V\right)H^{T}+I\right)^{-1}HG\right)^{-1}\le C.
\]
The last inequality stems from the fact that $P_{mix}\le \grave{P}_{\text{+}}$.

The block-diagonal approximation of the covariance of the intermediate filter is the solution of
\begin{equation*}
\tilde{\Pi}_{\text{+}}=\mathcal{ D}\left\{ \left(I-\tilde{\kappa}H\right)\left(F\tilde{\Pi}_{\text{+}}F^{T}+V+GUG^{T}\right)\left(I-\tilde{\kappa}H\right)^{T}+\tilde{\kappa}R\tilde{\kappa}^{T}\right\},
\end{equation*}
if it exists. It can be shown that the above equation is equivalent to ${\Pi}_{\text{+}}^{(i)}=\grave{P}_{\text{+}}^{(i)},\:i\ne1$ and 
\begin{equation*}
\Delta\tilde{\Pi}_{\text{+}}^{(1)}-\tilde{\Phi}_{c}^{(1)}\Delta\tilde{\Pi}_{\text{+}}^{(1)}\tilde{\Phi}_{c}^{(1)T}=F_{c}^{(1)}G^{(1)}C_{mix,1}G^{(1)T}F_{c}^{(1)T},
\end{equation*}
where $\Delta\tilde{\Pi}_{\text{+}}^{(1)}=\tilde{\Pi}_{\text{+}}^{(1)}-P_{\text{+}}^{(1)}$ and  $\tilde{\Phi}_{c}^{(1)}=\left(I-K_{mix,1}^{(1,1)}H^{(1)}\right)F^{(1)}$. This is a Lyapunov equation which has a solution if $\tilde{\Phi}_{c}^{(1)}$ is stable. To see that this is the case we note that
\begin{alignat*}{1}
\left\Vert \tilde{\Phi}_{c}^{(1)}-F_{c}^{(1)}\right\Vert _{F}= & \left\Vert \left(I+\left(F^{(1)}P_{\text{+}}^{(1)}F^{(1)T}+Q^{(1)}\right)H^{(1)T}\left(R^{(1)}\right)^{-1}H^{(1)}\right)^{-1}\cdot\right.\\
 & \cdot G^{(1)}C_{mix}G^{(1)T}\cdot\\
 & \left.\cdot H^{(1)T}\left(H^{(1)}\left(F^{(1)}P_{\text{+}}^{(1)}F^{(1)T}+V^{(1)}\right)H^{(1)T}+R^{(1)}\right)^{-1}H^{(1)}F^{(1)}\right\Vert _{F}\\
\le & \alpha_{3}^{(1)}\alpha_{4}^{(1)}\epsilon^{(1)}
\end{alignat*}
and that the condition $\max\left|\lambda\left\{ F_{c}^{(1)}(V)\right\} \right|+\beta^{(1)}\alpha_{3}^{(1)}\alpha_{4}^{(1)}\epsilon^{(1)}<1$ ensures that the maximum eigenvalue of $\tilde{\Phi}_{c}^{(1)}$ is still within the unit circle.

We will now find a bound on $\left\Vert \Delta\tilde{\Pi}_{\text{+}}^{(1)}\right\Vert $.
Let $\mathcal{F}_{X}$ be the discrete Lyapunov operator defined by
\begin{equation*}
\mathcal{F}_{X}^{(1)}\left\{ X\right\} =X-F_{c}^{(1)}XF_{c}^{(1)T}
\end{equation*}
and let 
\begin{equation}
\Delta P_{mix}^{(1)}=\mathcal{F}_{X}^{-1}\left\{ F_{c}^{(1)}G^{(1)}C_{mix}G^{(1)T}F_{c}^{(1)T}\right\}\label{eq:delta_P_mix}.
\end{equation}
Then, from a perturbation analysis of the Lyapunov equation\cite{gahinet1990sensitivity},
it follows that
\begin{equation}
\frac{\left\Vert \Delta\tilde{\Pi}_{\text{+}}^{(1)}-\Delta P_{mix}^{(1)}\right\Vert _{F}}{\left\Vert \Delta\tilde{\Pi}_{\text{+}}^{(1)}\right\Vert _{F}}\le\left(\left\Vert \left(\mathcal{F}_{X}^{(1)}\right)^{-1}\right\Vert \left(2\left\Vert F_{c}^{(1)}\right\Vert _{F}+\left\Vert \tilde{\Phi}_{c}^{(1)}-F_{c}^{(1)}\right\Vert _{F}\right)\right)\left\Vert \tilde{\Phi}_{c}^{(1)}-F_{c}^{(1)}\right\Vert _{F},\label{eq:gahinet}
\end{equation}
where $\left\Vert \left(\mathcal{F}_{X}^{(1)}\right)^{-1}\right\Vert $
is the norm induced by Frobenius norm and is given by,
\begin{equation*}
\left\Vert \left(\mathcal{F}_{X}^{(1)}\right)^{-1}\right\Vert =\left\Vert \left(\mathcal{I}-F_{c}^{(1)}\otimes F_{c}^{(1)}\right)^{-1}\right\Vert _{2}=\frac{1}{1-\max\left|\lambda\left\{ F_{c}^{(1)}\right\} \right|^{2}}=\alpha_{1}^{(1)},
\end{equation*}
$\mathcal{I}\in\mathbb{R}^{c^{2}\times c^{2}}$ is the identity matrix
and $\otimes$ is the Kroneker product. Finally, from eq.~(\ref{eq:delta_P_mix}) we have
\begin{equation*}
\left\Vert \Delta P_{mix,1}^{(1)}\right\Vert _{F}\le  \left\Vert \left(\mathcal{F}_{X}^{(1)}\right)^{-1}\right\Vert \left\Vert F_{c}^{(1)}\right\Vert _{2}^{2}\left\Vert G^{(1)}C_{mix,1}G^{(1)T}\right\Vert _{F}=\alpha_{1}^{(1)}\left(\alpha_{2}^{(2)}\right)^{2}\epsilon^{(1)}
\end{equation*}
and therefore if $\alpha_{1}^{(1)}\left(2\alpha_{2}^{(2)}+\alpha_{3}^{(1)}\alpha_{4}^{(1)}\epsilon^{(1)}\right)\alpha_{3}^{(1)}\alpha_{4}^{(1)}\epsilon^{(1)}<1$,
the left hand side of (\ref{eq:gahinet}) is also smaller than one and 
\begin{equation*}
\left\Vert \Delta\tilde{P}_{\text{+}}^{(1)}\right\Vert _{F}\le\left\Vert \Delta\tilde{\Pi}_{\text{+}}^{(1)}\right\Vert _{F}\le\frac{\left\Vert \Delta P_{mix}^{(1)}\right\Vert _{F}}{1-\frac{\left\Vert \Delta\tilde{\Pi}_{\text{+}}^{(1)}-\Delta P_{mix}^{(1)}\right\Vert _{F}}{\left\Vert \Delta\tilde{\Pi}_{\text{+}}^{(1)}\right\Vert _{F}}}=\frac{\alpha_{1}^{(1)}\left(\alpha_{2}^{(2)}\right)^{2}\epsilon^{(1)}}{1-\alpha_{1}^{(1)}\left(2\alpha_{2}^{(2)}+\alpha_{3}^{(1)}\alpha_{4}^{(1)}\epsilon^{(1)}\right)\alpha_{3}^{(1)}\alpha_{4}^{(1)}\epsilon^{(1)}}.
\end{equation*}

\subsubsection*{Proof of part (2)}

We consider the perturbation of the block diagonal approximation, eq.~(\ref{eq:diagonal_ricatti}), by following the analysis in \cite{konstantinov1993perturbation} (which applies to the full
KF as well). First we define the operator
\begin{equation*}
\mathcal{F}\left\{ X,Q\right\} =X-F\mathcal{D}\left\{ X-XH^{T}\left(HXH^{T}+R\right)^{-1}HX\right\} F^{T}-Q,
\end{equation*}
(note that $\mathcal{F}\left\{ \tilde{P}_{\text{- }},V+GUG^{T}\right\} =\mathcal{F}\left\{ P_{\text{-}},V\right\} =0$).
Its Frechet derivative w.r.t $X$ at $Q=V$ is
\begin{equation*}
\mathcal{F}_{X}\left\{ X\right\} =X-\mathcal{D}\left\{ F_{c}XF_{c}^{T}\right\} .
\end{equation*}
$\mathcal{F}_{X}$ is a linear operator over the space of matrices and
it can be ``described'' as
\begin{equation*}
\mathrm{vec}\left[\mathcal{F}_{X}\left\{ X\right\} \right]=\left(\mathcal{I}_{n^{2}}-\mathcal{D}_{n^{2}}\cdot\left(F_{c}\otimes F_{c}\right)\right)\cdot\mathrm{vec}[X]
\end{equation*}
where $\mathrm{vec}[\cdot]$ denotes stacking the columns of a matrix into a single vector, $\mathcal{I}_{n^{2}}\in\mathbb{R}^{(cn)^{2}\times(cn)^{2}}$
is an identitity matrix and $\mathcal{D}_{n^{2}}$ is a diagonal matrix
with only $0$ or $1$ entries. Since $F_{c}$ is stable, 
\begin{equation*}
\max\left|\lambda\left(\mathcal{D}_{n^{2}}\cdot\left(F_{c}\otimes F_{c}\right)\right)\right|\le\max\left|\lambda\left(F_{c}\otimes F_{c}\right)\right|=\max\left|\lambda\left(F_{c}\right)\right|^{2}<1,
\end{equation*}
hence $\mathcal{F}_{X}$ is invertible and
\begin{equation*}
\left\Vert \mathcal{F}_{X}^{-1}\left\{ Z\right\} \right\Vert _{F}\le\frac{1}{1-\max\left|\lambda\left\{ F_{c}\right\} \right|^{2}}\left\Vert Z\right\Vert _{F}=\alpha_{1}\left\Vert Z\right\Vert _{F},\:\forall Z.
\end{equation*}

Defining $\Delta P_{\text{-}}=\tilde{P}_{\text{-}}-P_{\text{-}}$,
one can show that
\begin{alignat*}{1}
0= & GUG^{T}+\mathcal{F}_{X}\left\{ \Delta P_{\text{-}}\right\} +\mathcal{B}\left\{ \Delta P_{\text{-}}\right\} ,\\
\mathcal{B}\left\{ Z\right\} = & \mathcal{D}\left\{ F_{c}ZH^{T}\left(H\left(P_{\text{-}}+Z\right)H^{T}+R\right)^{-1}HZF_{c}^{T}\right\} .
\end{alignat*}
The above equation can be rearranged to depict $\Delta P_{\text{-}}$ as a fixed point of $\mathcal{G}$, 
\begin{alignat*}{1}
\Delta P_{\text{-}}= & \mathcal{G}\left\{ \Delta P_{\text{-}}\right\}, \\
\mathcal{G}\left\{ Z\right\} = & -\mathcal{F}_{X}^{-1}\left\{ \mathcal{B}(Z)+GUG^{T}\right\},
\end{alignat*}
and we will use this fact to prove the bound on $\left\Vert \Delta P_{\text{-}}\right\Vert _{F}$.

From the definition of $\mathcal{B}$,
\begin{equation*}
\left\Vert \mathcal{B}\left\{ Z\right\} \right\Vert _{F}\le\left\Vert F_{c}\right\Vert _{2}^{2}\left\Vert H^{T}\left(HP_{\text{-}}H^{T}+R\right)^{-1}H\right\Vert _{2}\left\Vert Z\right\Vert _{F}^{2},\:\forall Z\ge0,
\end{equation*}
and therefore
\begin{equation*}
\left\Vert \mathcal{G}\left\{ Z\right\} \right\Vert _{F}\le\alpha_{1}\eta+\alpha_{1}\alpha_{2}^{2}\alpha_{5}\left\Vert Z\right\Vert _{F}^{2},\:\forall Z\ge0.
\end{equation*}
Thus, if $(\alpha_{1}\eta)\cdot(\alpha_{1}\alpha_{2}^{2}\alpha_{5})<\frac{1}{4}$,
$\mathcal{G}$ is a continuous function that maps the compact convex set $\left\{ \left.Z\ge0\right|\left\Vert Z\right\Vert _{F}^{2}\le\frac{1-\sqrt{1-4(\alpha_{1}\eta)\cdot(\alpha_{1}\alpha_{2}^{2}\alpha_{5})}}{2\alpha_{1}\alpha_{2}^{2}\alpha_{5}}\right\} $
to itself. It follows from Brouwer's fixed-point theorem that $\mathcal{G}$
has a fixed point in this set, which must be $\Delta P_{\text{-}}$ due to the uniqueness of
$\tilde{P}_{\text{-}}$. Furthermore, it can be shown
that
\begin{equation*}
F_{c}-\tilde{F}_{c}=\left(I+P_{\text{-}}H^{T}R^{-1}H\right)^{-1}\Delta P_{\text{-}}H^{T}\left(H\tilde{P}_{-}H^{T}+R\right)^{-1}HF,
\end{equation*}
and since $\tilde{P}_{\text{-}}\ge P_{\text{-}}$, we have 
\begin{equation*}
\left\Vert F_{c}-\tilde{F}_{c}\right\Vert \le\alpha_{3}\alpha_{4}\left\Vert \Delta P_{\text{-}}\right\Vert _{F}.
\end{equation*}
This proof for the full KF is almost identical.
\end{document}